%%%%%%%%%%%%%%%%%%%%%%%%%%%%%%%%%%%%%%%%%%%%%%%%%%%%%%%%%%%%%%%%%%%%%%%%%%%%
%%% On the curvature of a certain Riemannian space of matrices
%%% Attila Andai, Peter W. Michor, Denes Petz
%%% AmSTeX,  10 pages 23.09.99
%%%%%%%%%%%%%%%%%%%%%%%%%%%%%%%%%%%%%%%%%%%%%%%%%%%%%%%%%%%%%%%%%%%%%%%%%%%%
\input amstex
\input amsppt.sty
\hsize 35pc
\vsize 55pc
\def\nmb#1#2{#2}         % used for renumbering, TeX should ignore.
\def\cit#1#2{\ifx#1!\cite{#2}\else#2\fi} %for citing references
\def\totoc{}             %= to table of content, invoked by kms-book.sty
               % for producing index, invoked by kms-book.sty
\def\ign#1{}             %=ignore, invisible entry for the index only

\define\de{\delta}

\define\ze{\zeta}
\define\et{\eta}

\define\la{\lambda}
\define\lam{\lambda}

\define\Ga{\Gamma}

%%%%%%%%%%%% Denes' definition
\define\GG{G}

\define\F{F}
\def\<{\langle}
\def\>{\rangle}
\def\A{A}
\def\iK{{\Cal K}}
\def\bbbr{\Bbb R}

%%%%%%%%%%%% Attila's definition
\define\Tr{\operatorname{Tr}}

%%%%%%%%%%%%%%%%%%%%%%%%%
\redefine\i{^{-1}}
\define\x{\times}
\define\tr{\operatorname{Tr\,}}

\topmatter
\title
On the curvature of a certain Riemannian space of matrices
\endtitle
\author
Peter W. Michor, \\
D{\'e}nes Petz, \\ Attila Andai
\endauthor
\affil
Erwin Schr\"odinger International Institute of Mathematical Physics,
Wien, Austria
\endaffil
\address
P\.W\. Michor: Institut f\"ur Mathematik, Universit\"at Wien,
Strudlhofgasse 4, A-1090 Wien, Austria, and:
Erwin Schr\"odinger International Institute of Mathematical Physics,
Boltzmanngasse 9, A-1090 Wien, Austria
\endaddress
\email
peter.michor\@esi.ac.at
\endemail
\address
D\. Petz: Institute of Mathematics,
Technical University Budapest,
H-1521 Budapest XI. Szto\-czek u. 2, Hungary
\endaddress
\email
petz\@math.bme.hu
\endemail
\address
A\. Andai: Institute of Mathematics,
Technical University Budapest,
H-1521 Budapest XI. Szto\-czek u. 2, Hungary
\endaddress
\email
andaia\@math.bme.hu
\endemail
%\dedicatory \enddedicatory
\date {September 15, 1999} \enddate
%\thanks \endthanks
\abstract
The curvature tensor and the scalar curvature are computed in the
space of positive definite real matrices endowed by the Kubo-Mori
inner product as a Riemannian metric.
\endabstract
\endtopmatter

\document
\head Table of contents \endhead
%\input amspptb.sty
%\input \jobname.toc
%\loadtoc
%\loadindex
\block{\eightpoint
\noindent 1. Introduction \leaders \hbox to 
     1em{\hss .\hss }\hfill {\eightrm 1}\par 
\noindent 2. The Kubo-Mori metric on the space of positive definite 
     matrices \leaders \hbox to  1em{\hss .\hss }\hfill {\eightrm 1}\par 
\noindent 3. The submanifold of normalized matrices \leaders \hbox to 
     1em{\hss .\hss }\hfill {\eightrm 4}\par 
\noindent 4. Computation of the scalar curvature \leaders \hbox to 
     1em{\hss .\hss }\hfill {\eightrm 5}\par         
}\endblock

\head\totoc\nmb0{1}. Introduction \endhead
The state space of a finite quantum system is identified with the set of
positive semidefinite matrices of trace 1. The set of all strictly positive 
definite matrices of trace 1 becomes naturally a differentiable manifold
and the Kubo-Mori scalar product defines a Riemannian structure on it.
Reference \cit!{4} tells about the relation of this metric to the von Neumann 
entropy functional.

The objective of the paper is to compute the scalar curvature in the
Riemannian geometry of the Kubo-Mori scalar product. Actually, we consider
the space of real density matrices which is a geodetic submamifold in the
space of complex density matrices. Our study is strongly motivated by the
conjectures formulated in \cit!{3} and \cit!{4}. 
It was conjectured that the scalar
curvature takes its maximum when all eigenvalues of the density matrix
are equal, and more generally the scalar curvature is monotone with respect
to the majorization relation of matrices. Although we obtain an explicite 
formula for the scalar curvature, the conjecture remains unproven. 
Nevertheless, a huge number of numerical examples are still supporting 
the conjecture. The method of computation of this paper is inspired 
by \cit!{2}. 

When this paper was nearly finished we recieved the preprint \cit!{1} 
where the scalar curvature is computed for the Kubo-Mori
metric in the complex case by a different method.

\head\totoc\nmb0{2}. The Kubo-Mori metric on the space of positive
definite matrices \endhead

\subhead\nmb.{2.1}. The setup \endsubhead
Let $\Cal S=\Cal S(n)$ be the space of all real
selfadjoint $(n\x n)$-matrices, $\Cal S_+=\Cal S_+(n)$ be the open
subspace of positive
definite matrices. Then $\Cal S_+$ is a manifold with tangent bundle
$T\Cal S_+= \Cal S_+\x \Cal S$. We shall consider the following
Riemannian metric on $\Cal S_+$, where $D\in\Cal S_+$ and
$X,Y\in T_D\Cal S_+ = \Cal S$:
$$
G_D(X,Y) = \int_0^\infty \tr\big( (D+t)\i X(D+t)\i Y\big)\,dt.
$$
Because
$$
\Big|\tr\big( (D+t)\i X(D+t)\i Y\big) \Big| \le n t^{-2}\Vert X\Vert
\,\Vert Y \Vert
$$
the integral is finite. We shall identify $\Cal S$ with its dual $\Cal S^*$ by
the standard (i.e., Hilbert-Schmidt) inner product $\langle X,Y\rangle=
\tr(XY)$. Then we can view the Riemannian metric $G$ also as a mapping
$$\gather
\GG_D:T_D\Cal S_+ =\Cal S \to T^*_D\Cal S_+ = \Cal S^* \cong \Cal S,\\
\GG_D(X) = \int_0^\infty (D+t)\i X(D+t)\i\,dt,\\
\endgather$$
which is symmetric with respect to $\<X,Y\>= \tr XY$.
(Note that $\GG_D$ is the Frechet derivative of $\log D$.)
Now let $D\in \Cal S_+ $ and choose a basis of
$\Bbb R^n$ such that $D=\sum_i \la_i E_{ii}$ is diagonal, where
$(E_{ij})$ is the usual system of matrix units,
then the selfadjoint matrices
$$
\F_{kl} \equiv E_{kl}+E_{lk} \qquad (k \le l)
$$
are a complete system of eigenvectors of $\GG_D:\Cal S\to \Cal
S$. This means that $(F_{ij})_{1\le i\le j\le n}$ is an orthogonal basis
of $(T_D\Cal S, G_D)$ with
$$
G_D(F_{ij},F_{kl})=\cases 0 & \text{ for }(i,j)\ne (k,l) \\
                            2 m_{ij} &\text{ for }i=k<j=l\\
                            4 m_{ii} &\text{ for }i=j=k=l
\endcases
$$
where
$$
\int_0^\infty (\lambda_k+t)^{-1}  (\lambda_l+t)^{-1}\, dt=
\frac{\log \lambda_l - \log \lambda_k}{\lambda_l - \lambda_k}=:
{m_{kl}}\, .
$$
The expression $m_{kl}$ is a symmetric function of the eigenvalues
$\la_k$ and $\la_l$. In fact $1/m_{kl}$ is the logarithmic
mean of $\lam_k$ and $\lam_l$. This implies that
$$
m_{kl}=\frac{1}{\lam_k} \qquad \text{whenever} \qquad
\lam_k=\lam_l,
$$
in particular, $m_{kk}=1/\lam_k$. Note that $G_D(F_{ij})=m_{ij}F_{ij}$ 
for all $i\le j$.

Another symmetric expression
$$
\int_0^\infty (\lambda_i+t)^{-1}  (\lambda_j+t)^{-1}
(\lambda_k+t)^{-1}\, dt=\frac{m_{jk}-m_{ij}}{\lambda_i-\lambda_k}
=: {m_{ijk}}
$$
will appear below. The identity
$$
\frac{1}{m_{kl}}\left(\frac{m_{kkl}}{m_{kk}}+\frac{m_{kll}}{m_{ll}}\right)
=1
$$
is easily computed and will be used later.

\subhead\nmb.{2.2}. The Christoffel symbol \endsubhead
Since we have a global chart we can express the Levi-Civita
connection by one Christoffel symbol:
$$
(\nabla_\xi\et)|_D = d\et(D).\xi(D) - \Ga_D(\xi(D),\et(D))
$$
where $\xi,\et:\Cal S_+\to \Cal S$ are smooth vector fields.
The Christoffel symbol is then given by
$$
G_D(\Ga_D(X,Y),Z) = \tfrac12 dG(D)(Z)(X,Y)
     - \tfrac12 dG(D)(X)(Z,Y) - \tfrac12 dG(D)(Y)(X,Z),
$$
where the derivative of the metric
$$\multline
dG(D)(Z)(X,Y) =\int_0^\infty \tr\Bigl(-(D+t)\i Z(D+t)\i X(D+t)\i Y -\\
     - (D+t)\i X(D+t)\i Z(D+t)\i Y\Bigr)\,dt
\endmultline$$
is visibly symmetric in the entries $Z,X,Y$.

The Christoffel form is given by
$$\align
G_D(\Ga_D(X,Y)) &= -\tfrac12 dG(D)(X)(Y) \\
     &= \tfrac12 \int_0^\infty (D+t)\i
     \Bigl(X(D+t)\i Y+Y(D+t)\i X\Bigr)(D+t)\i\,dt,\\
\Ga_D(X,\quad) &= -\tfrac12 G_D\i.dG(D)(X) \\ .
%               &= \tfrac12 d(G\i)(D)(X).G_D,
\endalign$$
Since
$$
G\i_D(X)=\int_0^1 D^u X D^{1-u}\, du\, ,
$$
we can express the Christoffel form as an integral formula.
The derivative is
$$\align
d\Ga(D)(X)(Y,Z) &= -\tfrac12 d(G\i)(D)(X).dG(D)(Y)(Z)
     + \tfrac12 G_D\i\;d^2G(D)(X,Y)(Z).
\endalign$$

\subhead\nmb.{2.3} \endsubhead
When $D=\sum_i \lambda_i E_{ii}$ is diagonal, then
$$
dG(D)(Z)(X) = -\sum_{ijk} m_{ijk}\big( E_{ii}ZE_{jj}XE_{kk}+
E_{ii}X E_{jj} Z E_{kk}\big).
$$
In particular,
$$
dG(D)(F_{ij})(F_{kl}) =-\delta_{jk}m_{ilj}F_{il}-
\delta_{jl}m_{ikj}F_{ik}-
\delta_{il}m_{jki}F_{jk}-
\delta_{ik}m_{jli}F_{jl},
$$
and
$$
\Ga_D(F_{ij},F_{kl}) = -\delta_{jk}\frac{m_{ilj}}{m_{il}}F_{il}-
\delta_{jl}\frac{m_{ikj}}{m_{ik}}F_{ik}-
\delta_{il}\frac{m_{jki}}{m_{jk}}F_{jk}-
\delta_{ik}\frac{m_{jli}}{m_{jl}}F_{jl}
$$

\subhead\nmb.{2.4}. The curvature \endsubhead
The Riemannian curvature
$R(\xi,\et)\ze = (\nabla_\xi\nabla_\et -
\nabla_\et\nabla\xi -\nabla_{[\xi,\et]})\ze$ is then determined in
terms of the Christoffel form by
$$\align
R_D(X,Y)Z &= -d\Ga(D)(X)(Y,Z) + d\Ga(D)(Y)(X,Z) +\\
&\quad + \Ga_D(X,\Ga_D(Y,Z)) - \Ga_D(Y,\Ga_D(X,Z)).
\endalign$$
If we insert the expressions from \nmb!{2.2} we get after some
computation
$$\align
R_D(X,Y)Z &= \tfrac14 d(G\i)(D)(X).dG(D)(Y)(Z)
     -  \tfrac14 d(G\i)(D)(Y).dG(D)(X)(Z)\\
&= - \tfrac14 G_D\i.dG(D)(X).G_D\i.dG(D)(Y)(Z)\\
&\quad + \tfrac14 G_D\i.dG(D)(Y).G_D\i.dG(D)(X)(Z)
\endalign$$
The Ricci curvature is then given by the following trace
$$
\operatorname{Ric}_D(X,Z)
     = \tr_{\Cal S}\left(Y\mapsto R_D(X,Y)Z\right),
$$
and the scalar curvature is
$$
\operatorname{Scal}(D)
     = \tr_{\Cal S}(X\mapsto G_D\i.\operatorname{Ric}_D(X,\quad)).
$$

Next we compute the traces in a concrete basis.
Let $\A_s$ be an orthonormal basis with
respect to the inner product $\langle X,Y\rangle=\tr(XY)$ on $\Cal S$.
Then
$$
\operatorname{Ric}_D(X,Z)=\sum_s \< R_D(X,\A_s)Z,\A_s \>
$$
and
$$\align
\operatorname{Scal}(D)&=\sum_t \<G_D\i
.\operatorname{Ric}_D(\A_t,\quad),
\A_t\>= \sum_t \< \operatorname{Ric}_D(\A_t,\quad), G_D\i \A_t\>\\
&= \sum_t  \operatorname{Ric}_D(\A_t, G_D\i \A_t) = \sum_{t}\sum_s
\< R_D(\A_t,\A_s)G_D\i\A_t,\A_s\>\,.
\endalign$$

\head\totoc\nmb0{3}. The submanifold of normalized matrices  \endhead

\subhead\nmb.{3.1}. The submanifold of trace 1 matrices \endsubhead
We consider the affine submanifold of $\Cal S_+$ of all positive
definite real selfadjoint matrices with trace 1 and its tangent
bundle:
$$\align
\Cal S_1 &=\{D\in \Cal S_+: \tr(D)=1\},\\
T\Cal S_1 &= \Cal S_1\x \Cal S_0, \text{ where }
\Cal S_0 =\{X\in \Cal S: \tr(D)=0\}.
\endalign$$

\proclaim{Lemma}
\roster
\item
The unit normal field $n$ along the submanifold
$\Cal S_1$ with respect to the Riemannian metric $G$ from \nmb!{2.1}
is given by $n(D)=(D,D)$.
\item
The $G$-orthonormal projection
$P_D:\Cal S=T_D\Cal S_+ \to T_D\Cal S_1 = \Cal S_0$
is given by $P_D(X)=X-\tr(X)n(D)=X-\tr(X)D$ for $D\in \Cal S_1$ and
$X\in \Cal S_)$.
\item
The Christoffel form for the covariant derivative $\nabla^1$ of the induced
Riemannian metric $\Cal S_1, G^1$ is given by
$$
\Ga^1_D(X,Y) = P_D\Ga_D(X,Y) = \Ga_D(X,Y) -
\tr(\Ga_D(X,Y)).D,\quad D\in \Cal S_1,\; X,Y\in \Cal S_0,
$$
and the second fundamental form is given by
$$\align
S_D:{}& T_D\Cal S_1 \x T_D\Cal S_1 = \Cal S_0\x \Cal S_0 \to \Bbb R \\
S_D(X,Y) &= \tr(\Ga_D(X,Y)) \\
&= \int_0^\infty \tr\bigl( (D+t)\i X(D+t)\i Y \\
&\qquad\qquad  - \tfrac12 D(D+t)^{-2}X(D+t)\i Y\\
&\qquad\qquad  - \tfrac12 D(D+t)^{-2}Y(D+t)\i X\Bigr)\,dt.
\endalign$$
\endroster
\endproclaim

\demo{Proof}
If $X\in \Cal S$ commutes with $D$ we get
$$\align
G_D(X,Y) &= \int_0^\infty \tr((D+t)\i X(D+t)\i Y)dt =
     \int_0^\infty \tr((D+t)^{-2}XY)dt =\\
&= [-\tr((D+t)\i XY)]_{t=0}^{t=\infty} =\tr(D\i XY)
\endalign$$
Thus for $Y\in \Cal S_0$ we have
$G_D(D,Y) = \tr(Y) = 0$.
Moreover for $D\in \Cal S_1$ we have
$G_D(D,D) = \tr(D) = 1$, so \therosteritem1 follows.
The remaining assertions are standard facts from Riemannian geometry.

For the explicit expression of the second fundamental form we preceed
as follows. 
For $D\in \Cal S_1$ the Weingarten mapping is given by
$$\align
L_D:&T_D\Cal S_1 =\Cal S_0\to \Cal S_0=T_D\Cal S_1, \\
L_D(X) :&= \nabla_{(D,X)}n = dn(D).X - \Ga_D(X,n(D)) = X - \Ga_D(X,D),
\endalign$$
and the second fundamental form is then given by
$$\align
S_D:&T_D\Cal S_1 \x T_D\Cal S_1 = \Cal S_0\x \Cal S_0 \to \Bbb R \\
S_D(X,Y) :&= G_D(L_D(X),Y) = G_D(X - \Ga_D(X,D),Y) \\
&= G_D(X,Y) + \tfrac12 dG(D)(X)(D,Y) \\
&= \int_0^\infty \tr\bigl( (D+t)\i X(D+t)\i Y \\
&\qquad\qquad  - \tfrac12 D(D+t)^{-2}X(D+t)\i Y
     - \tfrac12 D(D+t)^{-2}Y(D+t)\i X\Bigr)\,dt.
\endalign$$
Another formula for the second fundamental form is
$$\align
S_D(X,Y) :&= \tr(\Ga_d(X,Y)) = \tr(-\tfrac12G_D\i\,dG(X)(Y))\tag4\\
&= \tfrac12\int_o^1\int_0^\infty \tr\bigl(
D^u(D+t)\i X(D+t)\i Y(D+t)\i D^{1-u} + \\
&\qquad\qquad  + D^u(D+t)\i Y(D+t)\i X(D+t)\i D^{1-u}\Bigr)\,dt\,du.
\qed\endalign$$
\enddemo

\subhead\nmb.{3.2}. The curvature via the Gau\ss{} equation\endsubhead
The Gau\ss{} equation expresses the curvature $R$ of $\Cal S_+$
and the curvature $R^1$ of $\Cal S_1$ for $D\in \Cal S_1$  and
$X,Y,Z,U \in \Cal S_0$ by
$$
G_D(R(X,Y)Z,U) = G_D(R^1(X,Y)Z,U) + S_D(X,Z)S_D(Y,U) - S_D(Y,Z)S_D(X,U).
$$

The Ricci curvature of the submanifold $\Cal S_1$
is then given by the following trace
$$
\operatorname{Ric}^1_D(X,Z)
     = \tr_{\Cal S_0}\left(Y\mapsto R^1_D(X,Y)Z\right),
$$
and the scalar curvature is
$$
\operatorname{Scal}^1(D)
     = \tr_{\Cal S_0}(X\mapsto G_D\i.\operatorname{Ric}^1_D(X,\quad)).
$$

Next we compute the traces in a concrete basis in case of a diagonal
$D=\sum_i \la_i E_{ii}\in \Cal S_1$. Let $\A_s$ be an orthonormal basis
with respect to the inner product $G_D$ on $\Cal S_0$.
Then
$$
\operatorname{Ric}^1_D(X,Z)=\sum_s G_D( R^1_D(X,\A_s)Z,\A_s )
$$
and
$$\align
\operatorname{Scal}^1(D)
&=\sum_t G_D(G_D\i.\operatorname{Ric}^1_D(\A_t,\quad),\A_t)
     = \sum_t \< \operatorname{Ric}^1_D(\A_t,\quad), \A_t\>\\
&= \sum_t  \operatorname{Ric}^1_D(\A_t,\A_t)
     = \sum_{t}\sum_s G_D( R^1_D(\A_t,\A_s)\A_t,\A_s)\\
& =\sum_{t,s} \Big( G_D( R_D(\A_t,\A_s)\A_t,\A_s)
-S_D(\A_t,\A_t)S_D(\A_s,\A_s)-
S_D(\A_s,\A_t)S_D(\A_t,\A_s)\Big)
\endalign$$

\head\totoc\nmb0{4}. Computation of the scalar curvature \endhead

Our aim is to have an explicit formula for the scalar curvature
$\operatorname{Scal}^1(D)$ in terms of eigenvalues the of $D$, which
we assume to be a diagonal matrix. As in
the previous section, let $\A_s$ be an orthonormal basis
with respect to the inner product $G_D$ on $\Cal S_0$. We assume
that some of the basis elements are diagonal (like $D$) and the others are
normalized symmetrized matrix units.

\subhead\nmb.{4.1}. The first term \endsubhead
We decompose the sum
$$
\sum_{t,s}  G_D( R_D(\A_t,\A_s) \A_t,\A_s)
$$
into three subsums and we compute them separately.
First we consider the case when both $\A_t$ and $\A_s$
are offdiagonal, that is, they are in the form $F_{ij}/\sqrt{2m_{ij}}$.

{\bf Offdiagonal-offdiagonal.}
$$\align
&\sum_{} \frac{1}{4m_{ij}m_{kl}} G_D( R_D(F_{ij},F_{kl}) F_{ij},
F_{kl})
=\sum_{} \frac{1}{4m_{ij}} \< R_D(F_{ij},F_{kl}) F_{ij}, F_{kl}\> \\
&=- \sum_{}
\frac{1}{16m_{ij}}\<G_D\i.d(G)(D)(\F_{ij}).G_D\i.dG(D)
(F_{kl})(\F_{ij}),F_{kl} \> \\
&\phantom{=} +\sum_{ijkl}
\frac{1}{16 m_{ij}} \<G_D\i.d(G)(D).(F_{kl}).G_D\i.dG(D)
(\F_{ij})(\F_{ij}),F_{kl} \>,
\endalign$$
where summation is over $i<j$ and $k<l$. We continue with the first term
and calculate in an elementary way:
$$\align &
-\sum_{} \frac{1}{16 m_{ij}}\<G_D\i.d(G)(D)(\F_{ij}).G_D\i.dG(D)
(F_{kl})(\F_{ij}),F_{kl} \> \\ & \qquad
\frac{12}{16}\sum_{u<v<w} \frac{m_{uvw}^2}{m_{uv}m_{vw}m_{uw}}
+\frac{2}{16}\sum_{i<j} \frac{m_{iij}^2}{m_{ij}^2 m_{ii}}
+\frac{2}{16}\sum_{i<j} \frac{m_{ijj}^2}{m_{ij}^2 m_{jj}}
\endalign $$

For the second term we use
$dG(D)(\F_{ij})(\F_{ij})=-m_{iij}F_{ii}-m_{ijj}F_{jj}$ and get
$$\multline
+\sum_{i<j,k<l}
\frac{1}{16 m_{ij}} \<G_D\i.d(G)(D).(F_{kl}).G_D\i.dG(D)
(\F_{ij})(\F_{ij}),F_{kl} \>\\
=\frac12\sum_{u<v<w}
\frac{m_{uuv}m_{uuw}}{m_{uv}m_{uw}m_{uu}} +
\frac12\sum_{u<v<w}
\frac{m_{vvw}m_{uvv}}{m_{vw}m_{uv}m_{vv}} +
\frac12\sum_{u<v<w}
\frac{m_{uww}m_{vww}}{m_{uw}m_{vw}m_{ww}} \\
+\frac14\sum_{i<j}
\frac{m_{iij}^2}{m_{ij}^2m_{ii}} +
\frac14\sum_{i<j}
\frac{m_{ijj}^2}{m_{ij}^2m_{jj}}
\endmultline$$

{\bf Offdiagonal-diagonal.} Next we compute the sum
$$
\sum_{t,s}  G_D( R_D(\A_t,\A_s) \A_t,\A_s)
$$
when $\A_t=\sum_ia^t_iE_{ii}$ are diagonal, $G_D$-orthogonal to $D$, and
orthonormalized, and where the $A_{s}$ are still offdiagonal.
This means that
$$\gather
G_D(D,\A_t) = \sum_i m_{ii}\la_i a^t_i =\sum_i a^t_i=0
\quad\text{ and }\quad
G_D(\A_t,\A_{t'}) = \sum_i m_{ii}a^{t}_i a^{t'}_i
=\sum_i \frac{a^{t}_i a^{t'}_i}{\la_i}=\de_{t,t'}.
\endgather$$
We also have
$$\gather
dG(D)(F_{kl})(A_{t})=dG(D)(A_t)(F_{kl}) = -(m_{kkl}a^t_k+m_{kll}a^t_l)F_{kl},
\quad dG(D)(A_t)(A_{t'}) = -D^{-2} A_tA_{t'},\\
dG(D)(F_{kl})(F_{kl})=-\delta_{kl}(m_{llk}+m_{kkl})F_{kl}
-m_{kkl}F_{kk}-m_{llk}F_{ll}, \quad\text{ and }
G_D\i A_t=DA_t
\endgather $$
since $D$, $A_t$, and $A_{t'}$ commute.
We get
$$\align
&\sum_{t,k<l}\frac1{2m_{kl}}
     G_D( R_D(\A_t,F_{kl}) \A_t,F_{kl}) \\
&=-\sum_{t,k<l}\frac1{8m_{kl}}
     \<dG(D)(\A_t).G_D\i.dG(D)
     (F_{kl})(\A_t),F_{kl} \> \\
&\phantom{=}+ \sum_{t,k<l}\frac1{8m_{kl}}
     \<dG(D)(F_{kl}).G_D\i.dG(D)
     (\A_t)(\A_t),F_{kl} \>=\\
&=-\sum_{t,k<l}\frac1{8m_{kl}}
     \<G_D\i.dG(D)(F_{kl})(\A_t),
     d(G)(D)(\A_t)(F_{kl}) \> \\
&\phantom{=}+ \sum_{t,k<l}\frac1{8m_{kl}}
     \<G_D\i.dG(D)(\A_t)(\A_t),
     dG(D)(F_{kl})(F_{kl}) \>=\\
&=-\sum_{t,k<l}\frac1{8m_{kl}}
     \<\frac1{m_{kl}}(-m_{kkl}a^t_k-m_{kll}a^t_l)F_{kl},
     (-m_{kkl}a^t_k-m_{kll}a^t_l)F_{kl}\> \\
&\phantom{=}+ \sum_{t,k<l}\frac1{8m_{kl}}
     \<-\sum_p \frac1{\la_p}(a^t_p)^2 E_{pp},
     (-m_{kkl}F_{kk}-m_{llk}F_{ll}) \>=\\
&=-\sum_{t,k<l}\frac1{4m_{kl}^2}
     (m_{kkl}a^t_k+m_{kll}a^t_l)^2 \\
&\phantom{=}+ \sum_{t,k<l}\frac1{4m_{kl}}
     \left(m_{kkl}\frac1{\la_k}(a^t_k)^2
     +m_{llk}\frac1{\la_l}(a^t_l)^2\right)=:Q\\
\endalign$$

Denoting by $Q$ this seemingly basis dependent quantity we transform the sums in 
$Q$ as follows:
$$
\sum_{t,k<l}(\dots)=\frac{1}{2}\left(\sum_{t,k,l}(\dots)-
\sum_{t,k=l}(\dots)\right).
$$
Summing for $k=l$ indexes we obtain:
$$\align
&-\sum_{t,k=l}\frac{1}{4m_{kl}^{2}}(m_{kkl}a_{k}^{t}+m_{kll}a_{l}^{t})^{2}=
-\sum_{t,k}\frac{m_{kkk}^{2}}{m_{kk}^{2}}(a_{k}^{t})^{2}=
-\frac{1}{4}\sum_{t,k}\frac{1}{\lambda_{k}^{2}}(a_{k}^{t})^{2}, \\
&\sum_{t,k=l}\frac{1}{4m_{kl}}\left(m_{kkl}\frac{1}{\lambda_{k}}(a_{k}^{t})^{2}+
m_{llk}\frac{1}{\lambda_{l}}(a_{l}^{t})^{2}\right)=
\sum_{t,k}\frac{1}{4m_{kk}}\frac{2m_{kkk}}{\lambda_{k}}(a_{k}^{t})^{2}=
\frac{1}{4}\sum_{t,k}\frac{1}{\lambda_{k}^{2}}(a_{k}^{t})^{2}.
\endalign $$
The two terms turned out to be equal, so
$$
Q=\frac{1}{8}\left[-\sum_{t,k,l}\frac{1}{m_{kl}^{2}}
(m_{kkl}a_{k}^{t}+m_{kll}a_{l}^{t})^{2}+\sum_{t,k,l}\frac{1}{m_{kl}}
\left(\frac{m_{kkl}}{\lambda_{k}}(a_{k}^{t})^{2}+
\frac{m_{llk}}{\lambda_{l}}(a_{l}^{t})^{2}\right)\right].
$$

We start to deal with the first sum. Let
$$
b^t_k={a^t_k \over \sqrt{\lambda_k}}, \quad \Lambda_k=\sqrt{\lambda_k}
\qquad (1 \le k \le n).
$$
Then $\Lambda,b^1,b^2,\dots,b^{(n-1)}$ is an orthonormal basis in $\bbbr^n$.
We define a linear mapping $\iK$ from $\bbbr^n$ to the space of all real $n\times n$
matrices (endowed by the standard Hilbert-Schmidt inner product).
$$
\iK c=\sum_{k,l}{1 \over m_{kl}}\big[m_{kkl}\sqrt{\lambda_k}c_k+m_{kll}
\sqrt{\lambda_l}c_l\big]E_{kl}\qquad (c\in \bbbr^n).
$$
Then
$$
\Vert \iK b^t\Vert^2=\sum_{k,l}{1 \over m_{kl}^2}\big[m_{kkl}a^t_k
+m_{kll}a^t_l\big]^2
$$
which is a term in $Q$. Hence
$$
\sum_{t,k,l}\frac{1}{m_{kl}^{2}}(m_{kkl}a_{k}^{t}+m_{kll}a_{l}^{t})^{2}=
\sum_t \Vert \iK b^t\Vert^2=\Tr \iK^*\iK-\Vert \iK \Lambda\Vert^2.
$$
Since
$$
\Vert \iK \Lambda\Vert^2=\sum_{k,l}{[m_{kkl}\lambda_k+m_{kll}\lambda_l]^2 
\over m_{kl}^2}=\sum_{k,l}1=n^2
$$
and
$$\align
\Tr \iK^*\iK& =\sum_i \Vert \iK e_i\Vert^2=
\sum_{i,l}2{m_{iil}^2\lambda_i\over m_{il}^2}+2 \sum_i {m_{iii}^2\lambda_i
\over m_{ii}^2}
\\ &= 2\sum_{k,l}{m_{kkl}^2\lambda_k\over m_{kl}^2}+{1 \over 2}\sum_i \lambda_i^{-1}
\endalign$$
we have
$$
\sum_{t,k,l}\frac{1}{m_{kl}^{2}}(m_{kkl}a_{k}^{t}+m_{kll}a_{l}^{t})^{2}=
2\sum_{k,l}{m_{kkl}^2\lambda_k\over m_{kl}^2}+{1 \over 2}\sum_i \lambda_i^{-1}
-n^2.
$$
The other terms of $Q$ are similarly computed as traces.
$$\align
\sum_{t,k,l}\frac{1}{m_{kl}}
\left(\frac{m_{kkl}}{\lambda_{k}}(a_{k}^{t})^{2}+
\frac{m_{llk}}{\lambda_{l}}(a_{l}^{t})^{2}\right) &=
\sum_{k,l}{m_{kkl}+m_{kll} \over m_{kl}}-
\sum_{k,l}{m_{kkl}\lambda_k+m_{kll}\lambda_l \over m_{kl}}\\&=
\sum_{k,l}{m_{kkl}+m_{kll} \over m_{kl}}-n^2
\endalign
$$
Finally, we obtain a basis independent expression for $Q$:
$$
Q=-\sum_{k,l}{m_{kkl}^2\lambda_k\over 4 m_{kl}^2}-
{1 \over 16}\sum_i \lambda_i^{-1}+
\sum_{k,l}{m_{kkl}+m_{kll}\over 8 m_{kl}}
$$

{\bf Diagonal-offdiagonal.} This case is completely similar (in fact, 
symmetric) and yields the same $Q$. 

{\bf Diagonal-diagonal.} Now we compute the sum
$$
\sum_{t,s}  G_D( R_D(\A_t,\A_s) \A_t,\A_s)
$$
when $A_t$ and $A_s$ are both diagonal,
$G_D$-orthogonal to $D$, and orthonormalized.
We get
$$\align
&\sum_{t,s}
     G_D( R_D(\A_t,A_s) \A_t,A_s) \\
&=-\sum_{t,s}\frac14
     \<d(G)(D)(\A_t).G_D\i.dG(D)
     (A_s)(\A_t),A_s \> \\
&\phantom{=}+ \sum_{t,s}\frac14
     \<d(G)(D)(A_s).G_D\i.dG(D)
     (\A_t)(\A_t),A_s \>=\\
&=-\sum_{t,s}\frac14
     \<G_D\i.dG(D)(A_s)(\A_t),
     d(G)(D)(\A_t)(A_s) \> \\
&\phantom{=}+ \sum_{t,s}\frac14
     \<G_D\i.dG(D)(\A_t)(\A_t),
     d(G)(D)(A_s)(A_s) \>=\\
&=-\sum_{t,s}\frac14
     (\<-D\i A_tA_s,-D^{-2}A_tA_s\>
     -\<-D\i(A_t)^2,-D^{-2} (A_s)^2\>)=0
\endalign$$

\subhead\nmb.{4.2}. The second term \endsubhead
We start the computation of
$$
-\sum_{t,s}S_D(\A_t,\A_t)S_D(\A_s,\A_s).
$$
We use first the formula from \nmb!{3.1}.\thetag4, and use also
\nmb!{2.3}.
$$\align
S_D(F_{ij},F_{kl}) &= \tr(-\tfrac12G_D\i\,dG(F_{ij})(F_{kl}))\\
&= \tr(-\tfrac12G_D\i(-\delta_{jk}m_{ilj}F_{il}-
     \delta_{jl}m_{ikj}F_{ik}-\delta_{il}m_{jki}F_{jk}-
     \delta_{ik}m_{jli}F_{jl}))\\
&= \tfrac12\tr\left(\delta_{jk}\frac{m_{ilj}}{m_{il}}F_{il}+
     \delta_{jl}\frac{m_{ikj}}{m_{ik}}F_{ik}
     +\delta_{il}\frac{m_{jki}}{m_{jk}}F_{jk}+
     \delta_{ik}\frac{m_{jli}}{m_{jl}}F_{jl}\right)\\
&= \delta_{jk}\de_{il}\frac{m_{ilj}}{m_{il}}+
     \delta_{jl}\de_{ik}\frac{m_{ikj}}{m_{ik}}
     +\delta_{il}\de_{jk}\frac{m_{jki}}{m_{jk}}+
     \delta_{ik}\de_{jl}\frac{m_{jli}}{m_{jl}}\\
\endalign$$
We observe that for $(i<j)\ne(k<l)$ we get $S_D(F_{ij},F_{kl})=0$.
Furthermore, for $i<j$  we have
$$\align
S_D(F_{ij},F_{ij}) &= \frac{m_{iij}}{m_{ii}}+\frac{m_{ijj}}{m_{jj}}\\
S_D(F_{ii},F_{kk}) &= 4\de_{ik}\frac{m_{iki}}{m_{ik}}\\
S_D(F_{ii},F_{kl}) &= \delta_{ik}\de_{il}\frac{m_{ili}}{m_{il}}+
     \delta_{il}\de_{ik}\frac{m_{iki}}{m_{ik}}
     +\delta_{il}\de_{ik}\frac{m_{iki}}{m_{ik}}+
     \delta_{ik}\de_{il}\frac{m_{ili}}{m_{il}}\\
&= 0 \text{ if }k<l.
\endalign$$

First we take summation when both $A_t$ and $A_s$ are offdiagonal:
$$\align
&-\sum_{t,s}S_D(\A_t,\A_t)S_D(\A_s,\A_s)=\\
=&-\sum_{i<j,k<l}
  S_D\left(\frac1{\sqrt{2m_{ij}}}F_{ij},\frac1{\sqrt{2m_{ij}}}F_{ij}\right)
  S_D\left(\frac1{\sqrt{2m_{kl}}}F_{kl},\frac1{\sqrt{2m_{kl}}}F_{kl}\right)
\\=&-\sum_{i<j,k<l}\frac1{4m_{ij}m_{kl}}
     \left(\frac{m_{iij}}{m_{ii}}+\frac{m_{ijj}}{m_{jj}}\right)
     \left(\frac{m_{kkl}}{m_{kk}}+\frac{m_{kll}}{m_{ll}}\right)
\\=&-\frac{1}{4}\sum_{i<j,k<l}1=-\frac{n^2(n-1)^2}{16}.
\endalign
$$

Next we assume that $A_t$ is diagonal and $A_s$ is offdiagonal.
$$
S_D(\A_t,\A_t)=\frac{1}{2}\sum_{p}(a_p^t)^2 \lambda_p^{-1}=\frac{1}{2},
$$
and
$$
-\sum_{t,s}S_D(\A_t,\A_t)S_D(\A_s,\A_s)=-\frac{1}{2}(n-1)\sum_{k<l}\frac1{4m_{kl}}
     \left(\frac{m_{kkl}}{m_{kk}}+\frac{m_{kll}}{m_{ll}}\right)=
-\frac{n(n-1)^{2}}{16}.
$$
Note that this contribution is symmetric in $t$ and $s$ and we should take it
twice.

Now it is easily follows the value of the sum when both $A_t$ and $A_s$
are diagonal:
$$
-\sum_{t,s}S_D(\A_t,\A_t)S_D(\A_s,\A_s)=
-\sum_{t,s}\frac{1}{2}\cdot\frac{1}{2} =-\frac{1}{4}(n-1)^2
$$

\subhead\nmb.{4.3}. The third term \endsubhead
Now we have to deal with
$$
-\sum_{t,s}S_D(\A_t,\A_s)S_D(\A_s,\A_t)=
-\sum_{t,s}S_D(\A_t,\A_s)^2
$$
By the formulas from \nmb!{4.2} we have that it suffices to sum when both $A_t$
and $A_s$ are diagonal and both of them are offdiagonal. Hence
$$
-\sum_{t,s}S_D(\A_t,\A_s)S_D(\A_s,\A_t)=
-{n-1\over 4} - \sum_{k<l}\frac{1}{4m_{kl}^2}
     \left(\frac{m_{kkl}}{m_{kk}}+\frac{m_{kll}}{m_{ll}}\right)^2=
-\frac{n-1}{4}-\frac{n(n-1)}{8}.
$$

\subhead\nmb.{4.4}. The scalar curvature formula \endsubhead
$$\align
R={}&\quad\frac{3}{4}\sum_{u<v<w}\frac{m_{uvw}^{2}}{m_{uv}m_{vw}m_{wu}}+
  \frac{1}{2}\sum_{u<v<w}\frac{m_{vvw}m_{vvu}}{m_{vu}m_{vv}m_{vw}}\\ 
&+\frac{1}{2}\sum_{u<v<w}\frac{m_{wwu}m_{wwv}}{m_{wu}m_{wv}m_{ww}}+
   \frac{1}{2}\sum_{u<v<w}\frac{m_{uuv}m_{uuw}}{m_{uu}m_{uv}m_{uw}}\\
&+\frac{3}{8}\sum_{i<j}\frac{m_{iij}^{2}}{m_{ij}^{2}m_{ii}}+
    \frac{3}{8}\sum_{i<j}\frac{m_{ijj}^{2}}{m_{ij}^{2}m_{jj}}\\
&-\sum_{k,l}{m_{kkl}^2\lambda_k\over 2 m_{kl}^2}-
{1 \over 8}\sum_i \lambda_i^{-1} 
+\sum_{k,l}{m_{kkl}+m_{kll} \over 4 m_{kl}} \\   
&-{n^2(n-1)^2\over 16}-{n(n-1)^2\over 8}-{1\over 4} (n-1)^2-{n-1\over 4}-
{n^2(n-1)^2\over 16}
\endalign $$

Some further simplification:
$$
\frac{3}{8}\sum_{i<j}\frac{m_{iij}^{2}}{m_{ij}^{2}m_{ii}}+
\frac{3}{8}\sum_{i<j}\frac{m_{ijj}^{2}}{m_{ij}^{2}m_{jj}} =
\frac{3}{16}\sum_{i,j}\frac{m_{iij}^{2}\lambda_i+m_{ijj}^{2}\lambda_j}
{m_{ij}^2}-{3\over 2}\sum_i {1 \over \lambda_i }
$$

$$
-{1 \over 2}\sum_{k,l}{m_{kkl}^2\lambda_k\over m_{kl}^2}=
-{1 \over 4}\sum_{k,l}{m_{kkl}^2\lambda_k+ m_{kll}^2\lambda_l\over m_{kl}^2}
$$

Hence
$$\align
R={}&\quad\frac{3}{4}\sum_{u<v<w}\frac{m_{uvw}^{2}}{m_{uv}m_{vw}m_{wu}}\\
&+ \frac{1}{2}\sum_{u<v<w}\frac{m_{vvw}m_{vvu}}{m_{vu}m_{vv}m_{vw}}+
   \frac{1}{2}\sum_{u<v<w}\frac{m_{wwu}m_{wwv}}{m_{wu}m_{wv}m_{ww}}+
   \frac{1}{2}\sum_{u<v<w}\frac{m_{uuv}m_{uuw}}{m_{uu}m_{uv}m_{uw}}\\
&- \frac{1}{16}\sum_{k,l}{m_{kkl}^2\lambda_k+ 
   m_{kll}^2\lambda_l\over m_{kl}^2}-
   \frac{13}{2}\sum_i {1 \over \lambda_i }+
   \sum_{k,l}{m_{kkl}+m_{kll} \over 4 m_{kl}}
   +\frac{n(n-1)}4(n^2-n+1).
\endalign$$

{\bf Acknowledgement.} The second-named author thanks to the Erwin 
Schr\"odinger International Institute for Mathematical Physics and to 
the Hungarian grant OTKA F023447 for support.


\Refs
\ref
\key{\cit0{1}}
\by J. Dittmann
\paper On the curvature of monotone metrics and a conjecture 
concerning the Kubo-Mori metric 
\paperinfo Preprint 1999
\endref

\ref
\key{\cit0{2}}
\by O. Gil-Medrano, P. W. Michor
\paper The Riemannian manifold of all Riemannian metrics
\jour Quaterly J. Math. Oxford (2)
\vol 42
\yr 1991
\pages 183--202
\endref

\ref
\key{\cit0{3}}
\by F. Hiai, D. Petz, G.Toth
\paper
Curvature in the geometry of canonical correlation
\jour
Studia Sci. Math. Hungar.
\vol 32
\yr 1996
\pages 235--249
\endref

\ref
\key{\cit0{4}}
\by D. Petz
\paper Geometry of canonical correlation on the state space of a
quantum system
\jour J. Math. Phys.
\vol 35
\yr 1994
\pages 780--795
\endref



\endRefs
\enddocument